\documentclass[reqno,12pt]{amsart}
\usepackage{graphicx}
\usepackage{epsfig}
\usepackage{amssymb,amsthm,amsmath,amsfonts,amstext,mathrsfs,fancybox}
\usepackage{enumerate}
\usepackage{latexsym}
\usepackage{epsfig}
\usepackage{url}
\usepackage{fancyhdr}
\usepackage{color}
\usepackage[unicode]{hyperref}

\hypersetup{colorlinks=true,
linkcolor=blue,
citecolor=blue,
urlcolor=blue,
filecolor=blue,
bookmarksnumbered=true,
pdfstartview=FitH,
pdfhighlight=/N}

\newtheorem{prop}{Proposition}
\newtheorem{prob}{Problem}
\newtheorem{defin}{Definition}

\input xy
\xyoption{all}

\def\d{\,{\rm{d}}}
\newcommand{\m}{\mathbf}

\title[Open problems]
{Algebraic functions with Fermat property, eigenvalues of transfer operator and Riemann zeros, and other open problems}
\author[G. Alkauskas]{Giedrius Alkauskas}
\address{Vilnius University, Department of Mathematics and Informatics, Naugarduko 24, LT-03225 Vilnius, Lithuania}
\email{giedrius.alkauskas@mif.vu.lt}

\begin{document}
\begin{abstract}
In this note we list a number of open problems in the fields of number theory, combinatorics, and representation theory: algebraic functions with Fermat property; power product expansion of the generating function for the partition function; relation between the non-trivial Riemann zeros and eigenvalues of the transfer operator; functional equation related to norm forms; two problems from geometric combinatorics; a problem on the moments of the Minkowski question mark function; a question in representation theory; a problem on interpolating the moments of the Stern's diatomic sequence; an arithmetic properties of the binary composition function.
\end{abstract}
\pagestyle{fancy}
\fancyhead{}
\fancyhead[LE]{{\sc Open problems}}
\fancyhead[RO]{{\sc G. Alkauskas}}
\fancyhead[CE,CO]{\thepage}
\fancyfoot{}
\date{September 30, 2016}
\thanks{The research of the author was supported by the Research Council of Lithuania grant No. MIP-072/2015.}

\maketitle

My research in the next couple of years will be concentrated on projective superflows, transfer operator for the Gauss continued fraction map (also the Mayer-Ruelle operator), and the integration of the world of Minkowski question mark function with the world of modular and quasi-modular forms. Moreover, more than half of the time I have commitments in areas outside science, and this time will increase after these two years. So, here the list of a number of open problems of various levels of difficulty from various areas of mathematics is presented. Hopefully some of these problems will catch an attention of professionals, and I hope that colleagues can propose them to students in mathematics.     

\section{Algebraic functions with a Fermat property \textbf{(number theory)}}

Consider the following rational function
\begin{eqnarray*}
P(x)=\frac{x^2}{(1-2x)(1-x)}=\sum\limits_{n=2}^{\infty}s(n)x^{n}=
\sum\limits_{n=2}^{\infty}(2^{n-1}-1)x^{n}.
\end{eqnarray*}
Thus, we have $s(p)\equiv 0\text{ (mod }p)$ for $p>2$ prime.\\

Next, let
\begin{eqnarray*}
J(x)=\frac{1}{\sqrt{1-4x}}-\frac{2}{1-x}=\sum\limits_{n=0}^{\infty}s(n)x^{n}=\sum\limits_{n=0}^{\infty}\Big{(}\binom{2n}{n}-2\Big{)}x^{n}.
\end{eqnarray*}
This also gives $s(p)\equiv 0\text{ (mod }p)$ for $p\geq2$ prime.\\

Finally, as is implied by (\cite{alkauskas-modular}, Proposition 1), we have the same (for $p>3$) conclusion for a degree $3$ algebraic function $T(x)=1+x^2+O(x^3)$, which satisfies
\begin{eqnarray*}
(6x^5-3x^4+2x^3+3x^2-1)T^3+(x^5-x^4+x^3+2x^2-1)T^2+(x^3-x^2+1)T+1=0.
\end{eqnarray*}
Let $R(x)=\sum\limits_{n=0}^{\infty}s(n)x^{n}\in\mathbb{Z}[[x]]$ be an algebraic function over $\mathbb{Q}(x)$, unramified at $x=0$.  Suppose, $s(p)\equiv 0\text{ (mod }p)$ for all sufficiently large prime numbers $p$. We call such a function $R(x)$ \emph{an algebraic function with a Fermat property}. Thus, we formulate
\begin{prob}
\label{prob1}
 Characterize all algebraic functions with a Fermat property in general, and in any particular algebraic function field $\mathbb{Q}(x,U)$, where $U$ is unramified at $x=0$.
\end{prob}
We can multiply the function $R(x)$ by an integer to get the congruence valid for all primes. All algebraic functions with a Fermat property form an abelian group $\mathscr{F}$. Since the set of all algebraic functions over $\mathbb{Q}$ is countable, $\mathscr{F}$ is also countable. If $U(x)=\sum_{n=0}^{\infty}a(n)x^{n}\in\mathbb{Z}[[x]]$ is an algebraic function, then $U^{\partial}(x):=xU'(x)=\sum_{n=0}^{\infty}na(n)x^{n}\in\mathscr{F}$. Indeed, first it is obvious that $U^{\partial}\in\mathbb{Z}[[x]]$. And second, if $Q(Y,x)\in\mathbb{Z}[Y,x]$ and $Q(U,x)=0$, then
\begin{eqnarray*}
U'=-\frac{Q_{x}(U,x)}{Q_{Y}(U,x)}
\end{eqnarray*} 
belongs to the same algebraic function field $\mathbb{Q}(U)$. Let $D$ (from ``Differential") be the union of all such possible $U^{\partial}$. Then $D$ is a subgroup of $\mathscr{F}$. So is $\mathbb{Z}[x]$. Finally, let $U(x)=\sum_{n=0}^{\infty}a(n)x^{n}\in\mathbb{Z}[[x]]$ is again an algebraic function, unramified and without a pole at $x=0$, and let for an integer $M\geq 2$, $U^{(M)}(x)=U(x^{M})$. Let $P$ (from ``Power") be the group whose elements are 
\begin{eqnarray*}
\sum\limits_{j=1}^{s}U_{j}^{(M_{j})},\quad M_{j}\geq 2.
\end{eqnarray*}
Then $P$ is also a subgroup of $\mathscr{F}$. Of course, any two of the subgroups $D$, $\mathbb{Z}[x]$ and $P$ have a non-trivial pairwise intersection. Let also for any $U$, algebraic over $\mathbb{Q}(x)$ and unramified at $x=0$, $\mathscr{F}_{U}=\mathscr{F}\cap\mathbb{Q}(x,U)$, $D_{U}=D\cap\mathbb{Q}(x,U)$, $P_{U}=P\cap\mathbb{Q}(x,U)$ (We identify any algebraic function with its Laurent expansion at $x=0$). We may refine Problem \ref{prob1} as follows.
\begin{prob}
\label{propb2}
Find the structure of abelian groups
\begin{eqnarray*}
\mathscr{F}/(D+\mathbb{Z}[x]+P),\quad \mathfrak{A}_{U}=\mathscr{F}_{U}/\big{(}D_{U}+\mathbb{Z}[x]\big{)},\quad
\mathfrak{P}_{U}=\mathscr{F}_{U}/\big{(}D_{U}+\mathbb{Z}[x]+P_{U}\big{)}
\end{eqnarray*}
for any $U$ algebraic over $\mathbb{Q}(x)$ and unramified at $x=0$. In particular, for example, what is the group $\mathfrak{P}_{\varnothing}$ (that is, we talk only about $\mathbb{Q}(x)$)?
\end{prob}
In fact, for a function $J(x)$ we have an even stronger property $s(p)\equiv 0\text{ (mod }p^{2})$. We may call it \emph{an algebraic function with a strong Fermat (or Wieferich) property}, and ask similar questions.\\
\noindent\rule{16.4cm}{1.2pt}

\section{Power product expansion of the generating function for the partition function \textbf{(combinatorics, number theory)}}
This problem is extracted from \cite{alkauskas-monthly}. Let $p(n)$ stand for the classical partition function. It is known that
\begin{eqnarray*}
P(x)=\sum\limits_{n=0}^{\infty}p(n)x^{n}=\frac{1}{\prod\limits_{k=1}^{\infty}(1-x^{k})}.
\end{eqnarray*}
Consider now the power product expansion of $P(x)$. That is, define the unique sequence of integers $n_{k}$, $k\geq 1$, such that
\begin{eqnarray*}
\sum\limits_{n=0}^{\infty}p(n)x^{n}=\prod\limits_{k=1}^{\infty}(1+n_{k}x^{k}).
\end{eqnarray*} 
This is the sequence A220420 in the Online Encyclopedia of Integer Sequences \cite{oeis}. It is easily demonstrated \cite{alkauskas-monthly} that $n_{k}=1$ if $k$ is odd, and for even indices, the sequence reads as 
\begin{eqnarray*}
2,4,0,14,-4,-8,-16,196,-54,-92,-184,144,-628,-1040,-2160,41102,\ldots
\end{eqnarray*}
Looking at these terms, we may ask the following.
\begin{prob}Explore the arithmetic properties of the sequence $n_{k}$. In particular, is it true that
\begin{eqnarray*}
n(8k+2),n(8k+4),n(8k+6)\text{ are negative, }n(8k)\text{ are positive for }k\geq 1.
\end{eqnarray*}
\end{prob}
This is probably easy (if we use the recurrence given in \cite{alkauskas-monthly}), but deeper properties of the sequence $\{n_{k}:n\in\mathbb{N}\}$ wait to be discovered.\\
\noindent\rule{16.4cm}{1.2pt}
\section{Map from set of non-trivial zeros of the Riemann zeta function to the set of eigenvalues of the transfer operator for the Gauss continued fraction map \textbf{(analytic number theory)}}
Let $\mathbb{D}$ be the disc $\{z\in\mathbb{C}:|z-1|<\frac{3}{2}\}$. Let $\mathbf{V}$ be the Banach space of functions which are analytic in $\mathbb{D}$ and are continuous in its closure, with the supremum norm.
For a complex number $s$, $\Re(s)>\frac{1}{2}$, one defines \emph{the Mayer-Ruelle transfer operator} by \cite{mayer}
\begin{eqnarray*}
L_{s}[f(t)](z)=\sum\limits_{m=1}^{\infty}\frac{1}{(z+m)^{2s}}f\Big{(}\frac{1}{z+m}\Big{)}.
\end{eqnarray*}
The operator $L_{s}$ is extended to all complex numbers $s$ by an analytic continuation. As the fundamental contribution, it was proved by D. Mayer that \cite{mayer}
\begin{eqnarray*}
\det(1-L_{s}^{2})=\prod\limits_{n=1}^{\infty}\Big{(}1-\lambda^{2}_{n}(s)\Big{)}=Z(s),
\end{eqnarray*}
where on the left we have the Fredholm determinant of the operator defined as the product involving eigenvalues $\lambda_{n}(s)$ of the operator $L_{s}$ (as given in the middle), and on the right - the Selberg zeta function for the full modular group.\\

It is known that one subset of non-trivial zeros of $Z(s)$ are given by $\frac{\rho}{2}$, where $\rho$ are non-trivial zeros of the Riemann zeta function \cite{zagier}. It was proved in \cite{alkauskas-transfer} (the project is being planned to be finished by March 2017), that the labelling of eigenvalues $\lambda_{n}$ with an integer $n$ is canonical, corresponding to polynomials of degree $(n+1)$: this ordering corresponds exactly to ordering $\lambda_{n}$ according to their magnitude.  In particular, for each nontrivial zero $\rho$ of the Riemann zeta function there exist an integer $M=t(\rho)$ such that
\begin{eqnarray*}
\lambda_{M}\Big{(}\frac{\rho}{2}\Big{)}=1.
\end{eqnarray*}
See \cite{zagier} to see that this number is $1$, which corresponds to Riemann zeros and even Maass wave forms,  not $-1$, which corresponds to odd Maass wave forms. For example, if we order nontrivial zeros according to the magnitude of imaginary part, calculations show that
\begin{eqnarray*}
t(\rho_{1})=1,\quad t(\rho_{2})=2,\quad t(\rho_{3})=1,\quad t(\rho_{4})=3,\quad t(\rho_{5})=1,\quad t(\rho_{6})=3.
\end{eqnarray*}
In particular, we pose
\begin{prob} Given an integer $M\in\mathbb{N}$. What one can say about the set $t^{-1}(M)$? Is it infinite? How the conjectural distribution of Riemann zeros change if we limit ourselves to the set $t^{-1}(M)$? At least, how one can compute, for example, the set $t^{-1}(1)$ effectively?
\end{prob}     
The question about trivial zeros of $Z(s)$ seems also interesting. Let $k\in\mathbb{N}$. It is known that  order of vanishing at $s=1-k$ of $Z(s)$ equals the dimension of the corresponding space of cusp forms $M_{2k}$ for ${\sf PSL}_{2}(\mathbb{Z})$. In particular,
\begin{prob}
How the number $\dim_{\mathbb{C}}(M_{2k})$ distributes among different factors of 
\begin{eqnarray*}
Z(1-k)=\prod\limits_{n=1}^{\infty}\Big{(}1-\lambda^{2}_{n}(1-k)\Big{)}.
\end{eqnarray*}
\end{prob}  
We re-iterate that the crucial thing here that the index $n$ is attributed canonically.\\
\noindent\rule{16.4cm}{1.2pt}
\section{Pseudo-automorphisms of norm forms in number fields \textbf{(algebraic number theory, functional equations)}}
This problem is taken form \cite{alkauskas-norm-form}.
\begin{prob} Let $T(a_{1},a_{2},\ldots,a_{n})$ be a norm form in some
integral basis of some proper field extension of $\mathbb{Q}$ of degree $n$.
Find all functions $f:\mathbb{Z}\rightarrow\mathbb{C}$, such that
\begin{eqnarray*}
T(f(a_{1}),f(a_{2}),\ldots,f(a_{n})) \text{ depends only on the value of }
T(a_{1},a_{2},\ldots,a_{n}).
\end{eqnarray*}
\end{prob}
This problem was solved by U. Zannier in case we have a norm form $a^{2}+b^{2}$. In essence, there are $6$ independent solutions, given by a constant, a square function, and characters of orders $2,3,4,5$ (see \cite{alkauskas-norm-form}).\\
 
In particular, consider the polynomial $h(X)=X^{3}-3X+1$. Let $\alpha$ be any root of it. Then the splitting field of this polynomial is cubic, given by $\mathbf{K}=\mathbb{Q}(\alpha)$. Let $\alpha'=\alpha^2-2$. The set $\{1,\alpha,\alpha'\}$ is an integral basis of the ring of integers, and the norm form is
\begin{eqnarray*}
T(a,b,c)=\mathcal{N}(a+b\alpha+c\alpha')=a^{3}-b^{3}-c^{3}-
3ab^{2}-3ac^{2}+3abc+6b^{2}c-3bc^{2}.
\end{eqnarray*}  
\begin{prob}
Find all functions $f,h:\mathbb{Z}\rightarrow\mathbb{C}$, such that
\begin{eqnarray*}
T\big{(}f(a),f(b),f(c)\big{)}=h\big{(}T(a,b,c)\big{)}.
\end{eqnarray*}
\end{prob}  
In \cite{alkauskas-norm-form} it is shown that if $f\equiv h$, then $f(n)=n$, $f(n)=-n$, $f(n)\equiv 0$, $f(n)\equiv i$, or $f(n)\equiv -i$. But it is much more natural to pose a problem with $f\neq h$, since this definitely contains much more information on the arithmetic of the field. \\
\noindent\rule{16.4cm}{1.2pt}
%\section{$2$-dimensional Thue-Morse sequence}
%\noindent\rule{17cm}{1.2pt}
\section{Friendly paths (\textbf{geometric combinatorics)}}
Here I reproduce my problem from American Mathematical Monthly \cite{alkauskas-paths}, Problem 11484, which still has no solution.\\

An \emph{uphill} lattice path is the union of a (doubly infinite) sequence of directed line segments in $\mathbb{R}^{2}$, each connecting an integer pair $(a,b)$ to an adjacent pair either $(a,b+1), $ $(a+1,b)$. A \emph{downhill} lattice path is defined similarly, but with
$b-1$ in place of $b+1$, and a \emph{monotone} lattice path is an uphill or downhill lattice path.\\

\indent Given a finite set $P$ of points in $\mathbb{Z}$, a \emph{friendly} path is a monotone lattice path for which there are as many points in $P$ on one side of the path as on the other. (Points that lie on the path do not count.) 
\begin{prob}
\label{prob4}
\indent Is it true that for every odd-sized set of points there is a friendly path?
\end{prob}
\noindent\rule{16.4cm}{1.2pt}
\section{Lines crossing all squares of the board (\textbf{geometric combinatorics)}}
This is my problem from c. 1998. It is so natural, that possibly it might had occurred to other people as well.\\

Let $n\in\mathbb{Z}$. Given a square $n\times n$, divided into $n^{2}$ unit squares in the usual way. A line is said to \emph{intersect} a unit square, if they have a common point, which is an inner point of the unit square. 
\begin{prob}Find the minimal number of lines, such that each of $n^{2}$ unit squares is intersected. In particular, is it true that for $n\geq 3$ this number is equal to $n-1$?
\end{prob}
\noindent\rule{16.4cm}{1.2pt}
\section{Formula for the moments of the Minkowski question mark function (\textbf{complex dynamics, number theory)}}
This problem is taken from \cite{alkauskas-ramanujan}, where it is posed as a Conjecture with a positive answer.\\

The \emph{Minkowski question mark function} $?(x)$ is defined by
\begin{eqnarray*}
?([0,a_{1},a_{2},a_{3},\ldots])=2^{1-a_{1}}-2^{1-(a_{1}+a_{2})}+2^{1-(a_{1}+a_{2}+a_{3})}-\ldots,\quad x\in[0,1],
\end{eqnarray*}
where $x=[0,a_{1},a_{2},\ldots]$ stands for the representation of $x$ by a regular continued fraction. \\

Now, define the rational functions  $\mathbf{Q}_{n}(z)=\mathscr{B}_{n}(z)z^{-n-1}$, $n\geq 0$, (where $\mathscr{B}_{n}(z)$ are polynomials with rational coefficients, which are $p-$adic integers for $p\neq 2$) by
\begin{eqnarray*}
\mathbf{Q}_{0}(z)=-\frac{1}{2z},\text{ and recurrently by }
\mathbf{Q}_{n}(z)=\frac{1}{2}\sum\limits_{j=0}^{n-1}\frac{1}{j!}\cdot\frac{\d^{j}}{\d
z^{j}}\mathbf{Q}_{n-j-1}(-1)\cdot\Big{(}z^{j}-\frac{1}{z^{j+2}}\Big{)}.
\end{eqnarray*}
The sequence $\mathbf{Q}^{'}_{n}(-1)=\frac{\d}{\d z}\mathbf{Q}_{n}(z)|_{z=-1},n\in\mathbb{N}_{0}$, begins with $\frac{1}{2},-\frac{1}{2},1,-\frac{5}{2},\frac{25}{4},-16,43,-\frac{971}{8},\frac{1417}{4},\ldots$
\begin{prob}
Is the following true? The function $\Lambda(t)=\sum\limits_{n=0}^{\infty}\frac{\mathbf{Q}^{'}_{n}(-1)}{n!}t^{n}$ is an entire function, and 
\begin{eqnarray*}
\int_{0}^{\infty}\Lambda(t)e^{-t}\d t=\int\limits_{0}^{1}x^{2}\d ?(x).
\end{eqnarray*}
\end{prob}
The same can be asked about higher moments, which are encoded by the same rational functions $\mathbf{Q}_{n}$; see \cite{alkauskas-ramanujan}.\\ 
\noindent\rule{16.4cm}{1.2pt}
\section{Superflows with the smallest possible symmetry group (\textbf{group representations, algebra)}}
We give a definition of the superflow, but in fact the question itself is purely algebraic from the theory of representation of finite groups over $\mathbb{C}$.    
\begin{defin}
Let $n\in\mathbb{N}$, $n\geq 2$, and $\Gamma\hookrightarrow{\rm GL}(n,\mathbb{R})$ be an exact representation of a finite group, and we identify $\Gamma$ with the image. We call the flow $\phi(\m{x})$ \emph{the $\Gamma$-superflow}, if 
\begin{itemize}
\item[i)]there exists a vector field $\mathbf{Q}(\m{x})=(Q_{1},\ldots,Q_{n})\neq (0,\ldots,0)$ whose components are $2$-homogenic rational functions and which is exactly the vector field of the flow $\phi(\m{x})$, such that
\begin{eqnarray}
\gamma^{-1}\circ\mathbf{Q}\circ\gamma(\mathbf{x})=
\mathbf{Q}(\mathbf{x})\label{kappa}
\end{eqnarray}
 is satisfied for all $\gamma\in\Gamma$, and
\item[ii)] every other vector field $\mathbf{Q}'$ which satisfies (\ref{kappa}) for all $\gamma\in\Gamma$ is either a scalar multiple of $\mathbf{Q}$, or  its degree of a common denominator is higher than that of $\mathbf{Q}$. 
\end{itemize}
The superflow is said to be \emph{reducible or irreducible}, if the representation  $\Gamma\hookrightarrow{\rm GL}(n,\mathbb{R})$ (considered as a complex representation) is reducible or, respectively, irreducible.
\end{defin}
Let $n\geq 2$ be an integer, and let
\begin{eqnarray*}
\psi_{n}=\inf\limits_{\Gamma}\frac{|\Gamma|}{n!};
\end{eqnarray*}
the infimum is taken over all finite subgroups of $\mathrm{GL}(n,\mathbb{R})$ for which there exist a superflow. It was proved in \cite{alkauskas-super} that for $n\geq 3$, one has $\psi_{n}\leq 2$.
\begin{prob}
Is it true that for $n\geq 3$, one has $\psi_{n}=2$? 
\end{prob}
\noindent\rule{16.4cm}{1.2pt}
\section{Function whose moments are moments of the Stern's diatomic sequence (\textbf{number theory, calculus)}}
Let $a(n)$ stand form the Stern's diatomic sequence: $a(0)=0$, $a(1)=1$, and for $n\geq 1$, one has $a(2n)=a(n)$, $a(2n+1)=a(n)+a(n+1)$.\\

Let us define the ``moments"
\begin{eqnarray*}
Q^{(L)}(N)=\sum\limits_{n=2^{N}+1}^{2^{N+1}}a^{L}(n),\quad
L\in\mathbb{N}_{0}.
\end{eqnarray*}
It is known (and easily provable) that $Q^{(1)}(N)=\sum\limits_{n=2^{N}+1}^{2^{N+1}}a(n)=3^{N}$. Also \cite{physics}, 
\begin{eqnarray*}
Q^{(2)}(N)=\frac{1}{\sqrt{17}}\Big{[}\Big{(}\frac{5+\sqrt{17}}{2}\Big{)}^{N+1}-\Big{(}\frac{5-\sqrt{17}}{2}\Big{)}^{N+1}\Big{]}.
\label{prop1}
\end{eqnarray*}
\begin{prop}Sequence $Q^{(L)}(N)$, $N\in\mathbb{N}$, satisfies the linear recurrence
of degree $L$, and there exists a unique algebraic number
$\mu_{L}$ (half of an algebraic integer) and a constant $c_{L}$, such
that
\begin{eqnarray*}
\frac{Q^{(L)}(N)}{2^{N}}\sim c_{L}\mu^{N}_{L}.
\end{eqnarray*}
\label{prop2}
\end{prop}
So, from above, $\mu_{0}=1$, $\mu_{1}=\frac{3}{2}$,
$\mu_{2}=\frac{5+\sqrt{17}}{4}$. The Table \ref{table1} gives information on
the initial values of $\mu_{L}$.\small
\begin{table}
\begin{tabular}{|r | r | r | r|}
\hline
\multicolumn{4}{|c|}{\textbf{1. Sequence $\mu_{L}$}}\\
\hline
$L$ & $\mu_{L}$& Minimal polynomial of $2\mu_{L}$ & $G$\\
\hline
$0$& $1$&$\lambda-2$ & $Z_{1}$\\
$1$ & $3/2$&$\lambda-3$ & $Z_{1}$\\
$2$ & $\frac{1}{4}(5+\sqrt{17})$&$\lambda^{2}-5\lambda+2$ & $Z_{2}$\\
$3$ & $7/2$&$\lambda-7$ & $Z_{1}$\\
$4$ & $\frac{1}{4}(11+\sqrt{113})$&$\lambda^{2}-11\lambda+2$ & $Z_{2}$\\
$5$ & $\frac{1}{2}(7+4\sqrt{6})$&$\lambda^{2}-14\lambda-47$ & $Z_{2}$\\
$6$ & $\frac{1}{4}(10+\sqrt{265}+\sqrt{357+20\sqrt{265}})$&$\lambda^{4}-20\lambda^{3}-161\lambda^{2}-40\lambda+4$ & $D_{4}$\\
$7$ & $20.50916706...$&$\lambda^{3}-29\lambda^2-485\lambda-327$ & $S_{3}$\\
$8$ & $\frac{1}{4}(22+\sqrt{1801}+\sqrt{2277+44\sqrt{1801}})$&$\lambda^{4}-44\lambda^{3}-1313\lambda^{2}-88\lambda+4$ & $D_{4}$\\
$9$ & 50.69978074...& $\lambda^{3}-65\lambda^{2}-3653\lambda-3843$ & $S_{3}$\\
$10$ & $\frac{1}{4}(50+3\sqrt{1345}+\sqrt{14597+300\sqrt{1345}})$ & $\lambda^{4}-100\lambda^{3}-9601\lambda^{2}-200\lambda+4$ & $D_{4}$\\
$11$ & $126.5114484...$ & $\lambda^{4}-156\lambda^{3}-24882\lambda^{2}+83828\lambda+107529$ & $S_{4}$\\
$12$ & $200.4256707...$ & $\lambda^{6}-247\lambda^{5}-63659\lambda^{4}+797003\lambda^{3}...$ & $Z_{2}\wr S_{3}$\\
\hline
\end{tabular}
\caption{The moments of the Stern's diatomic sequence, minimal polynomials of doubled moments, and Galois groups. $Z$ - cyclic, $D$ - dihedral, $S$ - symmetric.}
\label{table1}
\end{table}
\normalsize
It is also known that the sequence $\mu^{1/L}_{L}$ is increasing, bounded, and
\begin{eqnarray*}
\lim_{L\rightarrow\infty}\mu^{1/L}_{L}=\frac{1+\sqrt{5}}{2}.
\end{eqnarray*}
\begin{prob}Does there exist a function $\mathbf{p}(x)$ (analytic? continuous of bounded variation? generalized - that is, a distribution?) such that
\begin{eqnarray*}
\int\limits_{0}^{\frac{1+\sqrt{5}}{2}}x^{L}\d\mathbf{p}(x)=\mu_{L},\quad L\geq 0. 
\end{eqnarray*}
\end{prob}
If it does, this function has algebraic moments, and encodes a lot of structural information about continued fractions.\\
\noindent\rule{16.4cm}{1.2pt}
\section{Binary composition function}
Let $\vartheta(n)$ denote the number of compositions (ordered partitions) of a positive integer into powers of $2$. This is the sequence A023359 in The Online Encyclopedia of Integer Seqences \cite{oeis}. The generating function satisfies
\begin{eqnarray*}
\sum\limits_{n=0}^{\infty}\vartheta(n)x^{n}=\frac{1}{1-\sum\limits_{k=0}^{\infty}x^{2^{k}}}.
\end{eqnarray*}
It was proved in \cite{alkauskas-compos} that for any integer $a$ (positive, zero or negative) there exist a limit in a $2$-adic topology
\begin{eqnarray*}
\lim\limits_{k\rightarrow\infty}\vartheta(2^{k}+a)=\Theta(a)\in\mathbb{Z}_{2},
\end{eqnarray*}
where the latter stands for the ring of integers of the $2$-adic number field $\mathbb{Q}_{2}$. For example,
\begin{eqnarray*}
\Theta(0)=2^{3}+2^{8}+2^{9}+\cdots.
\end{eqnarray*}
\begin{prob}
Describe arithmetic properties of numbers $\Theta(a)$. Is it rational or algebraic over $\mathbb{Q}$?
\end{prob}
\noindent\rule{16.4cm}{1.2pt}


\begin{thebibliography}{99}


\bibitem{alkauskas-ramanujan}{\sc G. Alkauskas}, Semi-regular continued fractions and an exact formula for the moments of the Minkowski question mark function, {\it Ramanujan J.} {\bf 25} (3) (2011), 359--367.

\bibitem{alkauskas-norm-form}{\sc G. Alkauskas}, Functional equation related to quadratic and norm forms, {\it Lithuanian Math. }J. {\bf 45}(2) (2005), 123--141.

\bibitem{alkauskas-transfer}{\sc G. Alkauskas}, Transfer operator for the Gauss' continued fraction map. I. Structure of the eigenvalues and trace formula, \url{http://arxiv.org/abs/1210.4083}.

\bibitem{alkauskas-monthly},{\sc G. Alkauskas}, A curious proof of Fermat's little theorem, {\it Amer. Math. Monthly} {\bf 116}(4) (2009), 362--364.

\bibitem{alkauskas-modular}{\sc G. Alkauskas}, The modular group and words in its two generators, \url{http://arxiv.org/abs/1512.02596}.

\bibitem{alkauskas-paths}{\sc G. Alkauskas}, Friendly Paths,
{\it Amer. Math. Monthly} {\bf 119}(2) (2012), 167--168.

\bibitem{alkauskas-super}{\sc G. Alkauskas}, Projective superflows. II.  $O(3)$ and the icosahedral group, \url{http://arxiv.org/abs/1606.05772}.

\bibitem{alkauskas-compos}{\sc G. Alkauskas}, Congruence properties of the function that counts compositions into powers of 2, {\it J. Integer Seq.} {\bf 13}(5) (2010), Article: 10.5.3, 9 p. 

\bibitem{zagier}{\sc J. Lewis, D. Zagier}, Period functions for Maass wave forms. I, {\it Ann. of Math. (2)} {\bf 153} (2001), no. 1, 191--258.

\bibitem{mayer}{\sc D. Mayer}, The thermodynamic formalism approach to Selberg's zeta function for ${\rm PSL}(2,\mathbb{Z})$, {\it  Bull. Amer. Math. Soc. (N.S.)} {\bf 25} (1) (1991), 55--60.

\bibitem{physics}{\sc M. Waldschmidt, P. Moussa, J. M. Luck and C. Itzykson (edts.)}, {\it From number theory to physics. Papers from the Meeting on Number Theory and Physics held in Les Houches, March 7–16, 1989}. Springer-Verlag, Berlin, 1992.

\bibitem{oeis}The Online Encyclopedia of Integer Sequences, Sequence A265434.



\end{thebibliography}
\end{document}